\documentclass[12pt]{article}
\usepackage{theorem,amsfonts}

\newtheorem{theorem}{Theorem}
\newtheorem{proposition}{Proposition}
\newtheorem{lemma}{Lemma}
\newtheorem{corollary}{Corollary}
\theorembodyfont{\upshape}
\newtheorem{definition}{Definition}

\newtheorem{proof}{Proof}

\newcommand{\bt}{\begin{theorem}}
\newcommand{\et}{\end{theorem}}
\newcommand{\bl}{\begin{lemma}}
\newcommand{\el}{\end{lemma}}
\newcommand{\bp}{\begin{proposition}}
\newcommand{\ep}{\end{proposition}}
\newcommand{\bd}{\begin{definition}}
\newcommand{\ed}{\end{definition}}
\newcommand{\bc}{\begin{corollary}}
\newcommand{\ec}{\end{corollary}}
\newcommand{\be}{\begin{enumerate}}
\newcommand{\ee}{\end{enumerate}}
\newcommand{\bo}{\begin{proof}}
\newcommand{\eo}{\end{proof}}

\title{Normed convergence property for hypergroups admitting an invariant
measure}

\author{C. R. E. Raja}

\date{}

\begin{document}
\maketitle

\let\epsi=\epsilon
\let\vepsi=\varepsilon
\let\lam=\lambda
\let\Lam=\Lambda 
\let\ap=\alpha
\let\vp=\varphi
\let\ra=\rightarrow
\let\Ra=\Rightarrow 
\let\Llra=\Longleftrightarrow
\let\Lla=\Longleftarrow
\let\lra=\longrightarrow
\let\Lra=\Longrightarrow
\let\ba=\beta
\let\ga=\gamma
\let\Ga=\Gamma
\let\un=\upsilon

\begin{abstract}
We prove that a hypergroup admitting a countable basis and an 
invariant Haar measure has normed convergence property if and only if it
is compact.
\end{abstract}



\bigskip

Let $K$ be a (locally compact) hypergroup and ${\cal P}(K)$ be the space
of all regular Borel probability measures on $K$ equipped with the weak
topology.  
For $x \in K$, $\delta _x$ denotes the probability measure concentrated at
$x$ and for $\mu \in {\cal M}^1(K)$, $x \mu$ and $\mu x$ denote 
$\delta _x*\mu$ and $\mu *\delta _x$ respectively.  Also, for $\mu _1 ,
\mu _2 \in {\cal P}(K)$, $\mu _1 \mu _2$ denotes $\mu _1 *\mu _2$.  

\bd
{\it A hypergroup $K$ is said to have normed convergence property if for
any
sequence $(\mu _n)$ in ${\cal P}(K)$, there exists a sequence $(x_n)$ in
$K$ such that $(x_n \nu _n)$ is relatively compact where $\nu _n = \mu _n
\cdots \mu _2 \mu _1$ for all $n \geq 1$. }
\ed

It is well-known that a locally compact group
has normed convergence property if and only if it is compact (see
Theorem 2.3.2 of \cite{H})
and a Moore hypergroup has normed convergence property if and only if it
is compact (see \cite{BH2}).  In this short note we extend the
afore-mentioned 
result to hypergroups admitting an invariant Haar measure.

Let $K$ be a hypergroup admitting a right-invariant Haar measure $m$.  
Let $L^p(K)$ be the space of p-integrable functions on $K$ for 
$1\leq p <\infty$.  
Let ${\cal M}^1(K)$ be the space of all non-negative measures on $K$ with
total measure less than are equal to $1$.  
Now, for any $\mu \in {\cal M}^1(K)$ and for $1\leq p\leq \infty$, 
define $P_\mu \colon L^p(K) \ra L^p(K)$ by
$$P_\mu (f) (x) = \int f(xy) d\mu (y)= \int f(y) d(x\mu )(y)$$ where
$f(xy) = \int f(z) d(\delta _x*\delta _y)(z)$, for all $x \in K$.  It is
easy to see that $P_\mu $ is a contraction on $L^p(K)$ for any 
$\mu \in {\cal M}^1(K)$ and any $1\leq p \leq \infty$.  The map $\mu
\mapsto P_\mu$ is an injection, so we may view ${\cal M}^1(K)$ as a
subspace of bounded linear operators on $L^2(K)$.

We now recall the following topologies on ${\cal M}^1(K)$:

\be
\item weak topology on ${\cal M}^1(K)$ is the weakest topology for which
the function $\mu \mapsto \mu (f)$ is continuous for all bounded
continuous functions $f$ on $K$;

\item vague topology on ${\cal M}^1(K)$ is the weakest topology for which
the function $\mu \mapsto \mu (f)$ is continuous for all continuous
functions $f$ with compact support;

\item weak operator topology (respectively, strong operator topology) on
${\cal M}^1(K)$ is the subspace topology inherited from the space of
bounded linear operators on $L^2(K)$ with weak operator topology
(respectively, strong operator topology).
\ee
 
For a locally compact group $G$, it is well-known that on 
${\cal M}^1(G)$ the vague and weak operator topologies coincide and on 
${\cal P}(G)$,
the weak, vague and weak operator and strong operator topologies coincide 
(see Lemma 6.1.23 of \cite{H}).  The following extension of this result
may be proved arguing as in Lemma 6.1.23 of \cite{H} and we omit the
details.

\bp\label{p1} 
Let $K$ be a hypergroup admitting a countable basis and a
right-invariant measure.  Then we have the following:
\be
\item on ${\cal M}^1(K)$, the vague and
weak operator topologies coincide;

\item on ${\cal P}(K)$, the weak,
vague and weak operator and strong operator topologies coincide. 
\ee
In particular, ${\cal M}^1(K)$ is compact in the weak operator topology. 
\ep

We next prove the following:

\bl\label{l1}
Let $K$ be as in Proposition \ref{p1}.  Suppose $(\mu _n)$ is a sequence
in ${\cal P}(K)$ such that $(x_n \mu _n )$ is relatively compact.  Then
there exists a $f \in L^2(K)$ such that $||P_{\mu _n}(f) ||_2 \not \ra 0$.  
\el

\bo
Since $(x_n \mu _n)$ is relatively compact, there exists a 
$\delta >0$ and a compact set $C$ such that $x_n \mu _n (C) >\delta $ for 
all $n \geq 1$.  

Let $f$ be a continuous function on $K$ with compact support such that 
$0\leq  f \leq 1$ and $f (x) = 1$ for all $x \in C$.  Then 
$$P_{\mu _n} (f) (x_n) = \int f(y ) d(x_n\mu _n )(y) >\delta$$ 
for $n \geq 1$.  Let $(h_j)$ be a compactly supported left approximate
identity in the space of continuous functions vanishing at infinity.  
Then as in Proposition 3.1 of \cite{DL}, we have 
$$||P_{\mu _n} (f)||_\infty \leq ||h_j*f-f||_\infty +
||P_{\mu _n}(f)||_2 ||h_j||_2$$ where $h_j *f (y) = \int h_j(x) f(xy)
dm(x)$.  Since $||P_{\mu _n}(f)|| _\infty >\delta >0$ for all $n \geq 1$, 
we get that $||P_{\mu _n} (f) ||_2 \not \ra 0$.
\eo

\bt\label{t1}
Let $K$ be a hypergroup admitting a right-invariant Haar measure and a
countable basis.  Then $K$ has normed convergence property if and only if
$K$ is compact. 
\et

\bo
It is well-known that compact hypergroups have normed convergence
property.  We now prove the 'only if' part.  Suppose $K$ has normed
convergence property.  Let $\mu \in {\cal P}(K)$ be a symmetric measure
such that the support of $\mu$ is $K$. Then
by assumption, there exists a sequence $(x_n)$ such that $(x_n \mu ^n)$ is
relatively compact.  By Lemma \ref{l1}, we get that 
$$||P_\mu ^n (f) ||_2 \not \ra 0$$ as $n \ra \infty$ for some
continuous function $f$ with compact support.  By
convergence of alternating sequences (see \cite{BC}), in the strong
operator
topology $$P_\mu ^{2n} \ra P$$ where $P$ is a bounded
linear operator on $L^2(K)$.  This implies that $P$ is a self-adjoint
projection.  By Proposition \ref{p1}, there exists a $\rho \in 
{\cal M}^1(K)$ such that $P=P_\rho$.  Since $P_\mu ^n (f) \not \ra 0$,
$\rho \not =0$ and $\rho ^2 = \rho$.  This means that $\rho$ is the
normalized Haar measure on a compact subhypergroup $H$ of $K$ 
(see \cite{BH1}).  
Since $P_\mu ^2 P=P$, we get that the support of $\mu ^2$
is contained in $H$ (see \cite{BH1}).  Since the support of $\mu $ is 
$K$, we get that $K=H$ and hence $K$ is compact.  
\eo

\noindent {C. Robinson Edward Raja, \\ 
Stat-Math Unit,\\ 
Indian Statistical Institute,\\ 
8th Mile Mysore Road,\\
R. V. College Post,\\
Bangalore - 560 059.\\
India.  
e-mail: creraja@isibang.ac.in \\
or \\
e-mail: raja$\_$robinson@hotmail.com}

\end{document}